\documentclass[11pt,letter]{article}
\usepackage{amsfonts,epsf,amsmath,amssymb,graphicx,tabularx,booktabs}

\newtheorem{theorem}{\bf Theorem}[section]
\newtheorem{corollary}[theorem]{\bf Corollary}
\newtheorem{lemma}[theorem]{\bf Lemma}

\newcommand{\proof}{\noindent{\bf Proof.\ }}
\newcommand{\qed}{\hfill $\blacksquare$ \bigskip}

\begin{document}

\title{\bf{The weighted vertex PI index}}

\author{
Aleksandar Ili\' c \\
Faculty of Sciences and Mathematics, University of Ni\v s, Serbia \\
e-mail: \tt{aleksandari@gmail.com} \\
\and
Nikola Milosavljevi\' c \\
Faculty of Sciences and Mathematics, University of Ni\v s, Serbia \\
e-mail: \tt{nikola5000@gmail.com} \\
}

\date{\today}
\maketitle 

\begin{abstract}
The vertex PI index is a distance--based molecular structure descriptor, that recently found
numerous chemical applications. In order to increase diversity of this topological index for
bipartite graphs, we introduce weighted version defined as $PI_w (G) = \sum_{e = uv \in E} (deg (u)
+ deg (v)) (n_u (e) + n_v (e))$, where $deg (u)$ denotes the vertex degree of $u$ and $n_u (e)$
denotes the number of vertices of $G$ whose distance to the vertex $u$ is smaller than the distance
to the vertex~$v$. We establish basic properties of $PI_w (G)$, and prove various lower and upper
bounds. In particular, the path $P_n$ has minimal, while the complete tripartite graph $K_{n/3,
n/3, n/3}$ has maximal weighed vertex $PI$ index among graphs with $n$ vertices. We also compute
exact expressions for the weighted vertex PI index of the Cartesian product of graphs. Finally we
present modifications of two inequalities and open new perspectives for the future research.
\end{abstract}

 \vskip 0.1cm

 \vskip 0.1cm

\section{Introduction}
\label{sec:intro}

Let $G = (V, E)$ be a connected simple graph with $n = |V|$ vertices and $m = |E|$ edges. For
vertices $u, v \in V$, the distance $d (u, v)$ is defined as the length of the shortest path
between $u$ and $v$ in $G$.  The maximum distance in the graph $G$ is its diameter,
denoted by $d$.

In theoretical chemistry molecular structure descriptors (also called topological indices) are used
for modeling physico-chemical, pharmacologic, toxicologic, biological and other properties of
chemical compounds \cite{GuPo86}. There exist several types of such indices, especially those based
on vertex and edge distances \cite{KhAzAsWa09}. Arguably the best known of these indices is the
Wiener index~$W$, defined as the sum of distances between all pairs of vertices of the molecular
graph~\cite{DoEnGu01},
$$
W (G) = \sum_{u, v \in V} d (u, v).
$$

Besides of use in chemistry, it was independently studied due to its relevance in social science,
architecture and graph theory. With considerable success in chemical graph theory, various
extensions and generalizations of the Wiener index are recently put forward.

One of the oldest degree-based graph invariants are the first and the second Zagreb
indices~\cite{GuDa04}, defined as follows
\begin{eqnarray*}
M_1(G) &=& \sum_{u \in V(G)} deg(v)^2 \\
M_2(G) &=& \sum_{uv \in E(G)} deg(u) deg(v).
\end{eqnarray*}
The Zagreb indices and their variants have been used to study molecular complexity,
chirality, in QSPR and QSAR analysis, etc.

Let $e = uv$ be an edge of the graph $G$. The number of vertices of $G$ whose distance to the
vertex $u$ is smaller than the distance to the vertex $v$ is denoted by $n_u (e)$. Analogously,
$n_v (e)$ is the number of vertices of $G$ whose distance to the vertex $v$ is smaller than the
distance to the vertex $u$. The vertex PI index and Szeged index are defined as follows:
\begin{eqnarray*}
PI (G) &=& \sum_{e \in E} n_u (e) + n_v (e) \qquad \cite{DaGu10,Ha10,Il10,NaFaAs09}\\
SZ (G) &=& \sum_{e \in E} n_u (e) \cdot n_v (e) \qquad \cite{ChWu09,DaGu09,HeTa09,Il10a,IlKlMi10}
\end{eqnarray*}

In order to increase diversity for bipartite graphs, we introduce weighted versions of $PI$ and
$SZ$ index. In this paper we establish some basic properties of the weighted vertex $PI$ index and
prove various lower and upper bounds. We also present modifications of IMO 1984 and IMO 1999
inequalities and use them for establishing a sharp upper bound of weighted PI index. In addition,
we compute exact expressions for the weighted vertex PI index of the Cartesian product of graphs
and open new perspectives for the future research.

\section{Weighted version of vertex PI index}
\label{sec:2}

Let $P_n$ and $S_n$ denote the path and the star on $n$ vertices, and let $K_{n, m}$ denote the
complete bipartite graph.

We define the weighted version of the vertex PI index as follows
$$
PI_w (G) = \sum_{e \in E} (deg (u) + deg (v)) (n_u (e) + n_v (e))
$$

For bipartite graphs it holds $n_u (e) + n_v (e) = n$, and therefore the diversity of the original
$PI$ and $SZ$ indices is not satisfying. The following inequality holds for a graph $G$ with $n$
vertices and $m$ edges \cite{KhAzAsWa09}
$$
PI (G) \leq n \cdot m,
$$
with equality if and only if $G$ is bipartite. This is why we introduced weighted version of these
indices. Assume that every edge $e = uv$ has weight $deg (v) + deg (u)$. Now, if $G$ is a bipartite
graph, we have
\begin{equation}
\label{eq:bipartite} PI_w (G) = n \sum_{v \in V} deg^2 (v).
\end{equation}

This means that the weighted vertex PI index is directly connected to the first Zagreb index.
Furthermore, it follows that among bipartite graphs path $P_n$ and complete bipartite graph
$K_{\lfloor n/2 \rfloor, \lceil n/2 \rceil}$ have the minimum and maximum value of weighted vertex
PI index, respectively~\cite{GuDa04}. These values are
\begin{eqnarray*}
PI_w (P_n) &=& n(4n - 6) \\
PI_w (K_{\lfloor n/2 \rfloor, \lceil n/2 \rceil}) &=& n^2 \lfloor n/2 \rfloor \lceil n/2 \rceil.
\end{eqnarray*}

Next we present a new formula for computing the weighted vertex PI index of a graph.

\begin{lemma}
\label{le:formula} Let $G$ be a connected graph. Then $PI_w (G) = \sum_{x \in V} w_x (G)$, where
$$
w_x = \sum_{e = uv \in E, \ d (x, v) \neq d (x, u)} deg (u) + deg (v).
$$
\end{lemma}

\proof We apply double counting to the set of ordered pairs $(x, e)$ for which $x \in V (G)$, $e =
uv \in E(G)$ and $d (x, v) \neq d(x, u)$. Let $n_e (G) = |\{x \in G \mid d (x, v) \neq d(x, u)\}|$.
Then by definition it follows $PI_w (G) = \sum_{e \in E} n_e (G) \cdot (deg (u) + deg (v))$. On the
other hand, we have $\sum_{e \in E} n_e (G) \cdot (deg (u) + deg (v)) = \sum_{x \in V} w_x (G)$ and
$PI_w (G) = \sum_{x \in V} w_x (G)$, as desired. \qed

\section{Lower bounds}
\label{sec:3}

\begin{theorem}
Let $G$ be a connected graph on $n$ vertices, $m$ edges and diameter $d$. Then,
$$
PI_w (G) \geq 4d^2 - 4d - 2 + 6 m,
$$
with equality if and only if $G \cong P_n$.
\end{theorem}

\proof For $n = 2$, inequality is obvious. Otherwise, for each edge $e = uv \in E(G)$, we have $n_u (e) + n_v (e) \geq 2$ and $deg (u) + deg (v)
\geq 3$. Let $P_{d + 1} = v_0 v_1 \ldots v_d$ be a diametrical path. Since the distance of the
vertices $v_i$ and $v_j$ in the graphs $G$ and induced subgraph $P_{d + 1}$ is equal to $|i - j|$,
we have
\begin{eqnarray*}
PI_w (G) &=& \sum_{e \in E(P_{d+1}) \cup E (G) \setminus E (P_{d+1})} (deg (u) + deg (v))
(n_u (e) + n_v (e)) \\
&\geq& PI_w (P_{d + 1}) + \sum_{e \in E (G) \setminus E (P_{d+1})} (deg (u) + deg (v))
(n_u (e) + n_v (e)) \\
&\geq& (d + 1)(4d - 2) + 2 \cdot 3 \cdot (m - d) \\
&=& 4d^2 - 4d - 2 + 6 m.
\end{eqnarray*}

The equality holds if and only if there are no other vertices than those from $P_{d + 1}$ and since
$G$ is connected -- it follows that $G \cong P_n$. \qed

\begin{theorem}
Let $G$ be a connected graph on $n$ vertices. Then,
$$
PI_w (G) \geq n (4n - 6),
$$
with equality if and only if $G \cong P_n$.
\end{theorem}

\proof Let $v$ be an arbitrary vertex from $G$. Denote with $d = \max \{ d (v, u) | u \in G \}$ the
eccentricity of $v$, and define layers
$$
L_i (v) = \{ u \in V (G) \mid d (v, u) = i \}, \qquad i = 0, 1, \ldots, d.
$$

This layer representation of a graph is the main idea of the breadth first search algorithm for
graph traversals \cite{CoLRS01}. The graph $G$ has exactly two types of edges: the edges between
vertices of $L_i$, $0 \leqslant i \leqslant d$, and the edges connecting the vertices of $L_i$ and
those of $L_{i + 1}$, $0 \leqslant i \leqslant d - 1$. Let $H$ denote the subgraph induced by the
edges of the second type. From Lemma \ref{le:formula}, it follows that each vertex $v$ contributes
to $PI_w (G)$ exactly for the sum of weights of the edges of the second type (we can ignore the
edges connecting the vertices from the same layer). Notice that $deg_H(u) \leq deg(u)$ and that $H$ is bipartite. Therefore
\begin{eqnarray*}
w_v &=& \sum_{uv \in E(H)} (deg(u) + deg(v)) \\
&\geq& \sum_{uv \in E(H)} (deg_H(u) + deg_H(v)) \\
&=& \sum_{v \in V} deg_H^2 (v) \geq 4n - 6,
\end{eqnarray*}
since $4n - 6$ is the minimum value of the first Zagreb index. Equality holds if and only if $G$
has no edges of the first type and $H$ is isomorphic to $P_n$. Finally, we have $PI_w (G) = \sum_{v
\in V}w_v \geq n(4n - 6)$ with equality if and only if $G \cong P_n$. \qed

\section{Upper bounds}
\label{sec:4}

Let $e = uv$ be an arbitrary edge, such that it belongs to exactly $t (e)$ triangles. In that case,
it easily follows
$$
n_u (e) + n_v (e) \leq n - t (e) \qquad \mbox{and} \qquad deg (u) + deg (v) \leq n + t (e).
$$
Therefore,
\begin{equation}
\label{eq:triangle} PI_w (G) \leq \sum_{e \in E} (n - t (e))(n + t (e)) = n^2 m - \sum_{e \in E}
t^2 (e).
\end{equation}

A complete multipartite graph $K_{n_1, n_2, \ldots, n_k}$ is a graph in which vertices are adjacent
if and only if they belong to different partite sets. Let $T_{n, r}$ be the Tur\'{a}n graph which
is a complete $r$--partite graph on $n$ vertices whose partite sets differ in size by at most one.
This famous graph appears in many extremal graph theory problems \cite{Tu41}. Nikiforov in
\cite{Ni11} established a lower bound on the minimum number of $r$-cliques in graphs with $n$
vertices and $m$ edges (for $r=3$ and $r=4$). Fisher in \cite{Fi89} determined sharp lower bound
for the number of triangles, while Razborov in \cite{Ra08} determined asymptotically the minimal
density of triangles in a graph of given edge density.

\begin{theorem}
\label{th:triangleBound} Let $G$ be a connected graph on $n$ vertices, $m$ edges and $t$ triangles.
Then,
\begin{eqnarray}
\label{eq:triangle1} PI_w (G) \leq n^2 m - \frac{9 t^2}{m},
\end{eqnarray}
with equality if and only if $G \cong K_{a, b}$ for $t = 0$, and $G \cong T_{n, r}$ for $r \mid n$
and $t > 0$.
\end{theorem}

\proof The inequality directly follows from \eqref{eq:triangle}, by applying Cauchy--Schwarz
inequality
$$
\sum_{e \in E} 1^2 \sum_{e \in E} t^2 (e) \geq \left ( \sum_{e \in E} t (e) \right)^2 = (3 t)^2,
$$
since every triangle is counted exactly 3 times. The equality holds if and only if $t (e) = t'$ for
every edge $e = uv \in E$.

Therefore, using \eqref{eq:triangle} it follows that the equality in \eqref{eq:triangle1} holds if
and only if $n_u (e) + n_v (e) = n - t'$ and $deg (u) + deg (v) = n + t'$ holds for all edges $e
\in E$.

For $t = 0$, we have $PI_w (G) = n^2 m$ if and only if $G$ is a complete bipartite graph $K_{a,
b}$, $a + b = n$. In order to prove this, let $e = uv$ be an arbitrary edge of $G$. Since $t = 0$
implies $deg (u) + deg (v) = n$ and there are no triangles, the neighbors of $u$ form one
independent vertex partition and the neighbors of $v$ form the other independent vertex partition
of a bipartite graph $K_{a, b}$. Again, using $deg (x) + deg (y) = n$ for an arbitrary edge $xy$,
it follows that each vertex of $G$ is adjacent to all vertices from other partition and therefore
$K_{a, b}$ is a complete bipartite graph.

Assume now that $t > 0$. Let $v$ be a vertex from $V (G)$ with neighbors $N (v) = \{v_1, v_2,
\ldots, v_s\}$. Each vertex $v_i$ is adjacent to all vertices that are not adjacent with $v$. Let
$U = V (G) \setminus N (v)$ and assume that there is an edge $e$ incident with the vertices from
$U$. In that case, the number of triangles for that edge $t (e)$ is greater than or equal to $s$,
which is impossible since $t(e) = t(vv_1) \leq s - 1$ . Therefore, the vertices from $U$ form an
independent set. Let $u \neq v$ be an arbitrary vertex from $U$. Since $G$ is connected, $u$ is
adjacent with some vertex $v_i \in N (v)$ and since $deg(u) + deg(v_i) = n + t' = deg(v) +
deg(v_i)$, it follows that $deg(u) = deg(v)$ and $N(u) = N(v)$. Finally, every vertex from $U$ is
adjacent with every vertex in $N(v)$.

The vertex $v$ was arbitrary chosen, and it follows that $G$ is isomorphic to a complete
multipartite graph $K_{n_1, n_2, \ldots, n_r}$. Since each edge belongs to exactly $t'$ triangles,
it easily follows that $n_1 = n_2 = \ldots = n_r = k$. For $r \mid n$, we have
$$
|E (T_{n, r})| = \left (1 - \frac{1}{r} \right) \frac{n^2}{2} \qquad \mbox{and} \qquad t (T_{n, r})
= \frac{n (n - n/r)(n - 2n/r)}{6}
$$
and
$$
PI_w (T_{n, r}) = \left (1 - \frac{1}{r} \right) \frac{n^2}{2} \cdot 2 \left (n - \frac{n}{r}
\right) \cdot 2 \frac{n}{r} = n^2 m - \frac{9 t^2}{m}.
$$

This completes the proof. \qed

Next we will establish a sharp upper bound for the weighed PI index. For that we need some
preliminary results.
\medskip

The following lemma is strongly connected to Problem 1 of International Olympiad in Mathematics
1984 \cite{DjJa06}.

\begin{lemma}
\label{lm:basic} Let $a, b, c$ be positive real numbers, such that $a + b + c = 1$. Then,
$$
0 \leq ab+bc+ac-abc \leq \frac{8}{27}.
$$
\end{lemma}

\proof From Cauchy--Schwarz inequality it follows
$$
(a + b + c)\left (\frac{1}{a} + \frac{1}{b} + \frac{1}{c} \right) \geq 9,
$$
and $ab+bc+ac \geq 9abc \geq abc$. For the right-side inequality, using AM--GM inequality we have
\begin{eqnarray*}
f (a, b, c) &=& ab+bc+ac-abc = c (a + b) + ab (1 - c) \leq c (1 - c) + \frac{(a + b)^2}{4} (1 - c)
\\&=& (1 - c) \cdot \frac{1 + c}{2} \cdot \frac{1+c}{2} \leq \left (\frac{1 - c + \frac{1 + c}{2} + \frac{1 + c}{2}}{3} \right)^3 = \left (\frac{2}{3} \right)^3 = \frac{8}{27}.
\end{eqnarray*}
with equality if and only if $a = b$ and $1 - c = \frac{1 + c}{2}$, i.e. $a = b
= c = \frac{1}{3}$. \qed

\begin{lemma}
\label{lm:squaresum} Let $a_1 \geq a_2 \geq \ldots \geq a_n \geq X > 0$ be positive real numbers,
such that $a_1 + a_2 + \ldots + a_n = Y$ and $n \geq 2$. Then
$$
\sum_{i = 1}^n{a_i}^2 \leq (Y - X)^2 + X^2,
$$
with equality if and only if $n = 2$ and $a_1 = Y - X$, $a_2 = X$.
\end{lemma}

\proof Notice first that $Y - (n - 1)X \geq X$. Suppose that $S(a) = \sum_{i = 1}^n{a_i}^2$ reaches
its maximum for some $n$-tuple $(a_1, a_2, \ldots, a_n) \neq (Y - (n-1)X, X, \ldots, X)$. Then,
there exist indices $i \neq j$ such that $a_j \geq a_i > X$. Let $\Delta = a_i - X$. By taking
$a_j' = a_j + \Delta$ and $a_i' = a_i - \Delta$, we increase the value of $S (a)$ since ${a_i'}^2 +
{a_j'}^2 - ({a_i}^2 + {a_j}^2) = 2 \Delta^2 + 2\Delta(a_j - a_i) > 0$. This is clearly a
contradiction.

Therefore for fixed $n$, the maximum is achieved for $a_1 = Y - (n-1)X$ and $a_i = X$, $i = 2, 3,
\ldots, n$, and its value is $S = (Y - (n-1)X)^2 + (n-1)X^2$. For $1 \leq k \leq n - 1$ let $S_k =
(Y - kX)^2 + kX^2$. Since $S_{k-1} - S_k = 2(Y - kX)X > 0$, we get
$$
S = S_{n-1} < S_{n-2} < \ldots < S_1 = (Y - X)^2 + X^2,
$$
which completes the proof. \qed

The following theorem is strongly connected to Problem 2 of International Olympiad in Mathematics
1999 \cite{DjJa06}.

\begin{theorem}
\label{th:multipart} Let $a_1 \geq a_2 \geq \ldots \geq a_n$ be positive real numbers, such that
$a_1 + a_2 + \ldots + a_n = 1$. Then,
$$
\sum_{i < j} a_i a_j (a_i + a_j) (2 - a_i - a_j) \leq \frac{8}{27},
$$
with equality if and only if $a_1 = a_2 = a_3 = \frac{1}{3}$ and $a_4 = \ldots = a_n = 0$.
\end{theorem}

\proof Let
$$
F (a_1, a_2, \ldots, a_n) = \sum_{i = 1}^{n - 1} \sum_{j = i + 1}^n a_i a_j (a_i + a_j) (2 - a_i -
a_j).
$$
First, we show the following inequality
\begin{equation}
\label{ieq:reduction}  F(a_1, a_2, \ldots, a_n) \leq F(a_1, a_2, \ldots, a_{n-1} + a_n, 0)
\end{equation}
for all $n \geq 4$. If $a_{n-1} = 0$ or $a_n$ = 0, the inequality is obvious. Otherwise, we have
\begin{eqnarray*}
\Delta &=& F(a_1, a_2, \ldots, a_n) - F(a_1, a_2, \ldots, a_{n-1} + a_n, 0) \\
&=& \sum_{i = 1}^{n-2}a_ia_n(a_i + a_n)(2 - a_i - a_n) + \sum_{i = 1}^{n-2}a_ia_{n-1}(a_i + a_{n-1})(2 - a_i - a_{n-1}) \\
& & + a_{n-1}a_n(a_{n-1} + a_n)(2 - a_{n-1} - a_n) \\
& & - \sum_{i = 1}^{n-2}a_i(a_{n-1} + a_n)(a_i + a_{n-1} + a_n)(2 - a_i - a_{n-1} - a_n) \\
&=& \sum_{i = 1}^{n-2}a_ia_{n-1}a_n(-4 + 4a_i + 3a_{n-1} + 3a_n) + a_{n-1}a_n(a_{n-1} + a_n)(2 - a_{n-1} - a_n) \\
&=& a_{n-1}a_n \left [ -4\sum_{i = 1}^{n-2}a_i + 4\sum_{i = 1}^{n-2}{a_i}^2 + 3(a_{n-1} + a_n)\sum_{i = 1}^{n-2}a_i + 2(a_{n-1} + a_n) - (a_{n-1} + a_n)^2 \right ] \\
\end{eqnarray*}
Let $x = a_{n-1}$ and $y = a_{n}$. Using $\sum_{i = 1}^{n-2}a_i = 1 - a_{n-1} - a_n = 1 - x - y$, we get
$$
\Delta = xy(-4 + 4\sum_{i = 1}^{n-2}{a_i}^2 + 9x + 9y - 4(x+y)^2).
$$
Therefore, $\Delta \leq 0$ is equivalent to
$$
\sum_{i = 1}^{n-2}{a_i}^2 \leq 1 + (x+y)^2 - \frac{9}{4}(x+y).
$$
Next, we will consider the following two cases.

{\bf Case 1.} $n = 4$. \\
Using $x + y = 1 - a_1 - a_2$ we get
$$
\sum_{i = 1}^{n-2}{a_i}^2 \leq 1 + (x+y)^2 - \frac{9}{4}(x+y) \quad \Leftrightarrow \quad a_1 + a_2 + 8a_1a_2 \geq 1.
$$
By ordering we have $a_2 \geq \frac{1 - a_1}{3}$ and after substitution it suffices to prove $(a_1
- 1)(4a_1 - 1) \leq 0$ which is true because $1 \geq a_1 \geq \frac{1}{4}$.

Equality occurs only for $(a_1, a_2, a_3, a_4) = (\frac{1}{4}, \frac{1}{4}, \frac{1}{4},
\frac{1}{4})$, since we exclude the case $(1, 0, 0, 0)$ because of the assumption $y \neq 0$.

{\bf Case 2.} $n > 4$. \\
By applying Lemma \ref{lm:squaresum} on $a_1, \ldots, a_{n-2}$ with $X = x$ and $Y = 1 - x - y$, we
get
$$
\sum_{i = 1}^{n-2}{a_i}^2 < (1 - 2x - y)^2 + x^2.
$$
Now, it suffices to prove
$$
(1 - 2x - y)^2 + x^2 \leq 1 + (x+y)^2 - \frac{9}{4}(x+y).
$$
Simplification gives $7x \geq 16x^2 + 8xy + y$, which can be easily verified using $x \geq y$ and
$x \leq \frac{1}{4}$. \medskip

With this, inequality \eqref{ieq:reduction} is proven. Notice that inequality is strict unless $a_n
= 0$, or $n = 4$ and $(a_1, a_2, a_3, a_4) = (\frac{1}{4}, \frac{1}{4}, \frac{1}{4}, \frac{1}{4})$.
Using this inequality and induction, we can reduce the problem to $n = 3$,
$$
F(a_1, a_2, a_3) = \sum_{1 \leq i < j \leq 3}{a_ia_j(1 - a_{6-i-j})(1 + a_{6-i-j})} = a_1a_2 + a_2a_3 + a_3a_1 - a_1a_2a_3.
$$
Now the result follows directly from Lemma \ref{lm:basic} and $F(\frac{1}{4}, \frac{1}{4},
\frac{1}{4}, \frac{1}{4}) = \frac{9}{32} < \frac{8}{27}$. \qed

The value of weighted PI index of complete multipartite graph can be easily calculated,
$$
PI_w (K_{n_1, n_2, \ldots, n_k}) = \sum_{i < j} n_i n_j (n_i + n_j) (2n - n_i - n_j).
$$
By substituting $a_i = \frac{n_i}{n}$ and applying Theorem \ref{th:multipart}, it follows that
among all multipartite graphs, the balanced 3-partite graph $K_{\lfloor n/3 \rfloor, \lceil n/3
\rceil, n - \lfloor n/3 \rfloor - \lceil n/3 \rceil}$ is the unique graph with maximum weighted
vertex PI index. Indeed, maximum is achieved for some partition of size 3. The case $3 \mid n$
follows directly from Theorem \ref{th:multipart}; otherwise assume that $(n_1, n_2, n_3)$ is some
3-partition with $n_3 - n_1 \geq 2$. It can be easily verified that
$$
PI_w (K_{n_1, n_2, n_3}) - PI_w (K_{n_1 + 1, n_2, n_3 - 1}) = n(n - n_2)(n_1 + 1 - n_3) < 0,
$$
which means that $(n_1, n_2, n_3)$ is balanced.
\medskip

For $3 \mid n$ it holds $PI_w(K_{n/3, n/3, n/3}) = \frac{8}{27} n^4$. We now show that
$\frac{8}{27} n^4$ is the upper bound for $PI_w(G)$ among all graphs on $n$ vertices and that this
bound is sharp for $3 \mid n$.

\begin{figure}[ht]
  \center
  \includegraphics [width = 6cm]{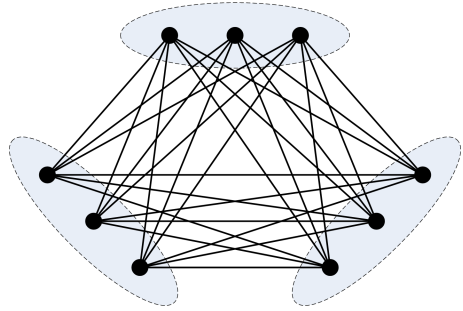}
  \caption { The complete tripartite graph $K_{3, 3, 3}$. }
\end{figure}

\begin{theorem}
Let $G$ be a connected graph on $n$ vertices. Then,
$$
PI_w (G) \leq \frac{8}{27} n^4,
$$
with equality if and only if $3 \mid n$ and $G \cong K_{n/3, n/3, n/3}$.
\end{theorem}

\proof Let $G_{n, m}$ be connected graph on $n$ vertices and $m$ edges. Denote by $t_{n, m}$ the
smallest possible number of triangles in $G_{n, m}$. From Theorem \ref{th:triangleBound} it follows
that
$$
PI_w(G_{n, m}) \leq n^2 m - \frac{9 t^2}{m} \leq n^2 m - \frac{9 t_{n,m}^2}{m}
$$
In \cite{NoSt63} the authors proved that
$$
t_{n,m} \geq \frac{(4m - n^2)m}{3n}
$$
with equality if and only if $G$ is a complete multipartite graph with partitions of equal size.
Substituting this in the previous inequality and after some simplification, we get
$$
PI_w(G_{n, m}) \leq 8 \cdot \frac{n^2m^2 - 2m^3}{n^2}.
$$
It suffices to prove that $8\cdot \frac{n^2m^2 - 2m^3}{n^2} \leq \frac{8}{27} n^4$, which is
equivalent to $n^6 + 54m^3 \geq 27n^2m^2$. This is true by simple AM--GM inequality:
$$
\left( \frac{n^2}{3} \right)^3 + m^3 + m^3 \geq n^2m^2,
$$
with equality if and only if $m = \frac{n^2}{3}$. It follows that $PI_w(G_{n, m}) \leq \frac{8}{27}
n^4$ for all $n - 1 \leq m \leq \binom{n}{2}$, and therefore $PI_w(G) \leq \frac{8}{27} n^4$ as
desired. From the above analysis, the equality holds if and only if $G$ is a regular multipartite
graph with partitions of size $\frac{n}{3}$. \qed

\section{Cartesian product graphs}

In this section we present formulas for computing weighted PI index of the Cartesian product of
graphs. In \cite{KlRaGu96, Kl07, KhAzAs08} the authors computed the PI index, the vertex PI index
and the Szeged index of Cartesian product graphs, respectively.

For graphs $G$ and $H$, the Cartesian product $G \times H$ is a graph with vertex set $V(G \times
H) = V(G) \times V(H)$ and $(u', u'')(v', v'')$ is an edge of $G \times H$ if $u' = v'$ and $u''v''
\in E(H)$, or $u'v' \in V(G)$ and $u'' = v''$. We will use the following well-known assertions for
the Cartesian product of graphs (see book of Imrich and Klav\v zar \cite{ImKl00} for more details)
and vertices $u = (u', u'')$, $v = (v', v'')$:
\begin{itemize}
\item $|V(G \times H)| = |V(G)| \cdot |V(H)|$
\item $|E(G \times H)| = |V(G)||E(H)| + |V(H)||E(G)|$
\item $deg_{G \times H} (u) = deg_G (u') + deg_H (u'')$
\item $d_{G \times H}(u, v) = d_G (u', v') + d_H (u'', v'')$.
\end{itemize}

Let us recall alternative formulas for computing the vertex and weighted PI index (see
Lemma~\ref{le:formula} and \cite{NaFaAs09}): $PI_v (G) = \sum_{x \in G} m_x(G)$ and $PI_w (G) =
\sum_{x \in G} w_x(G)$, where
\begin{eqnarray*}
w_x(G) &=& \sum_{e = uv \in E(G), \ d (x, v) \neq d (x, u)} deg (u) + deg (v), \\
m_x(G) &=& \sum_{e = uv \in E(G), \ d (x, v) \neq d (x, u)} 1.
\end{eqnarray*}

\begin{theorem}
\label{th:cartProd} Let $G$ and $H$ be two connected graphs. Then
$$
PI_w(G \times H) = |V(G)|^2 PI_w(H) + |V(H)|^2 PI_w(G) + 4 \Big( |V(G)||E(G)|PI_v(H) +
|V(H)||E(H)|PI_v(G) \Big).
$$
\end{theorem}

\proof We will use the formula
$$
PI_w(G \times H) = \sum_{x \in V(G \times H)} w_x(G \times H).
$$

Notice that for $x, u, v \in G \times H$, the condition $d(x, u) \neq d(x, v)$ is equivalent with
$$d((x', x''), (u', u'')) \neq d((x', x''), (v', v'')) \Leftrightarrow d(x', u') + d(x'', u'') \neq
d(x', v') + d(x'', v'').$$

If $(u, v) \in E(G \times H)$ it follows the either $u' = v'$ or $u'' = v''$, and therefore
\begin{eqnarray*}
w_x(G \times H) &=& \sum_{e = uv \in E(G \times H), \ d (x, v) \neq d (x, u)} deg (u) + deg (v) \\
&=& \sum_{u''v'' \in E(H), \ u' = v', \ d (x'', v'') \neq d (x'', u'')} deg(u') + deg(v') + deg(u'') + deg(v'') + \\
& & \sum_{u'v' \in E(G), \ u'' = v'', \ d (x', v') \neq d (x', u')} deg(u') + deg(v') + deg(u'') + deg(v'') \\
&=& A (x)+ B(x).
\end{eqnarray*}

We have
\begin{eqnarray*}
A (x) &=& \sum_{u' \in G} \phantom{a} \sum_{u''v'' \in E(H), \ d (x'', v'') \neq d (x'', u'')} 2 \cdot deg(u') + deg(u'') + deg(v'') \\
&=& \sum_{u' \in G} (m_{x''}(H) \cdot 2 \cdot deg(u') + w_{x''}(H)) \\
&=& 4|E(G)|m_{x''}(H) + |V(G)|w_{x''}(H).
\end{eqnarray*}

Analogously, $B (x) = 4|E(H)|m_{x'}(G) + |V(H)|w_{x'}(G)$ and thus
$$
w_x(G \times H) = |V(G)|w_{x''}(H) + |V(H)|w_{x'}(G) + 4|E(G)|m_{x''}(H) + 4|E(H)|m_{x'}(G).
$$

Finally, it follows
\begin{eqnarray*}
PI_w(G \times H) 
&=& \sum_{x' \in V(G), \ x'' \in V(H)} A(x) + B (x)\\
&=& |V(G)| \sum_{x'' \in V(H)} A (x) + |V(H)| \sum_{x' \in V(G)} B (x)\\
&=& |V(G)|(4|E(G)|PI_v(H) + |V(G)|PI_w(H)) \\
& & + |V(H)|(4|E(H)|PI_v(G) + |V(H)|PI_w(G)),
\end{eqnarray*}
which completes the proof.\qed

Denote by $\bigotimes_{i = 1}^n G_i$ the Cartesian product of graphs $G_1 \times G_2 \times \ldots
\times G_n$ and let $|V(G_i)| = V_i$ and $|E(G_i)| = E_i$ for all $i = \overline{1, n}$. Using the
above properties of Cartesian product graphs, one can easily verify that $V(\bigotimes_{i = 1}^n
G_i) = \prod_{i = 1}^n V_i$ and $E(\bigotimes_{i = 1}^n G_i) = \sum_{i = 1}^n E_i \prod_{j = 1, \ j
\neq i}^n V_i$. In \cite{KhAzAs08}, Khalifeh et al. have proven
\begin{equation}
\label{eq:cartProdVertex} PI_v(\bigotimes_{i = 1}^n G_i) = \sum_{i = 1}^n PI_v(G_i)\prod_{j = 1, \ j \neq i}^n |V(G_j)|^2.
\end{equation}
We prove similar result for weighted PI index:

\begin{theorem}
\label{th:cartProdGeneral} Let $G_1, G_2, \ldots, G_n$ be connected graphs. Then
$$
PI_w(\bigotimes_{i = 1}^n G_i) = \sum_{i = 1}^n PI_w(G_i)\prod_{j = 1, \ j \neq i}^n {V_j}^2 + 4\sum_{i, j = 1, \ i \neq j}^n PI_v(G_i)V_jE_j\prod_{k = 1, i \neq k \neq j}^n {V_k}^2.
$$
\end{theorem}

\proof The case $n = 2$ was proven in Theorem \ref{th:cartProd}. We continue our argument by
mathematical induction. Suppose that the result is valid for some $n$ graphs. Using Theorem
\ref{th:cartProd} and equation~\eqref{eq:cartProdVertex} we have
\begin{eqnarray*}
PI_w(\bigotimes_{i = 1}^{n+1} G_i) &=& PI_w(\bigotimes_{i = 1}^n G_i \times G_{n + 1}) \\
&=& |V(\bigotimes_{i = 1}^n G_i)|^2 PI_w(G_{n+1}) + V_{n+1}^2 PI_w(\bigotimes_{i = 1}^n G_i) \\
& & + 4 \Big (|V(\bigotimes_{i = 1}^n G_i)||E(\bigotimes_{i = 1}^n G_i)|PI_v(G_{n+1}) + V_{n+1}E_{n+1}PI_v(\bigotimes_{i = 1}^n G_i) \Big) \\
&=& PI_w(G_{n+1}) \prod_{i=1}^n {V_i}^2 + V_{n+1}^2 \sum_{i = 1}^n PI_w(G_i)\prod_{j = 1, \ j \neq i}^n V_i^2\\
& & + 4\sum_{i, j = 1, \ i \neq j}^n PI_v(G_i) V_j E_j \prod_{k = 1, i \neq k \neq j}^{n+1} V_k^2\\
& & + 4 \Big(PI_v(G_{n+1})\sum_{i=1}^n E_i V_i \prod_{j=1, \ j \neq i}^n V_j^2 + V_{n+1} E_{n+1} \sum_{i = 1}^n PI_v(G_i)\prod_{j = 1, \ j \neq i}^n V_j^2 \Big) \\
&=& \sum_{i = 1}^{n + 1} PI_w(G_i)\prod_{j = 1, \ j \neq i}^{n + 1} {V_j}^2 \\
& & + 4 \Big(\sum_{i, j = 1, \ i \neq j}^n PI_v(G_i) V_j E_j \prod_{k = 1, i \neq k \neq j}^{n+1} V_k^2 + \sum_{\substack{1 \leq i, j \leq n \\ i = n + 1 \vee j = n + 1}} PI_v(G_i) V_j E_j \prod_{k = 1, i \neq k \neq j}^{n+1} V_k^2 \Big) \\
&=& \sum_{i = 1}^{n + 1} PI_w(G_i)\prod_{j = 1, \ j \neq i}^{n + 1} {V_j}^2 + 4\sum_{i, j = 1, \ i
\neq j}^{n + 1} PI_v(G_i)V_jE_j\prod_{k = 1, i \neq k \neq j}^{n + 1} {V_k}^2.
\end{eqnarray*}
This completes the proof. \qed

\begin{corollary}
Let $G$ be connected graph. Then
$$
PI_w(G^n) = PI_w(\bigotimes_{i = 1}^n G) = n|V(G)|^{2n - 3} \Big(|V(G)|PI_w(G) + 4
(n-1)|E(G)|PI_v(G)\Big).
$$
\end{corollary}
\proof Directly follows from Theorem \ref{th:cartProdGeneral}. \qed

\section{Concluding remarks}

The vertex PI and Szeged indices are novel molecular-structure descriptors. In this paper we
generalized these indices and open new perspectives for the future research. Similarly we can
define weighted Szeged index as follows
$$
SZ_w (G) = \sum_{e \in E} (deg (u) + deg (v)) n_u (e) \cdot n_v (e).
$$

It would be interesting to study mathematical properties of these modified indices and report their
chemical relevance and formulas for some important graph classes. In particular, some exact
expressions for the weighted PI and SZ index of other graph operations (such as the composition,
join, disjunction and symmetric difference of graphs, bridge graphs, Kronecker product of graphs)
can be derived, similarly as in \cite{HoLuVu10,KhAzAs08a,KlRaGu96,MaSc09}.

\bigskip {\bf Acknowledgement. } This work was supported by Research
Grants 174010 and 174033 of Serbian Ministry of Science.

\end{document}